\documentclass[12pt]{amsart}
\input{epsf}

\usepackage{amssymb,latexsym}
\topskip  \numberwithin{equation}{section}
\newtheorem{theorem}{Theorem}[section]

\newtheorem{corollary}[theorem]{Corollary}
\newtheorem{definition}[theorem]{Definition}
\newtheorem{conjecture}[theorem]{Conjecture}

\newtheorem{lemma}[theorem]{Lemma}

\newtheorem{example}[theorem]{Example}

\newtheorem{problem}[theorem]{Problem}

\newcommand{\E}{\mathcal{E}}

\newcommand{\sym}{S}

\newcommand{\e}{\varepsilon}

\newcommand{\PR}{\mathcal{PR}}
\newcommand{\Au}{\mathcal{A}}

\newcommand{\ST}{\mathcal{ST}}
\def\ch#1,#2,{\binom#1#2}
\def\hy#1,#2,#3,#4,{{}_2F_1\left(#1,#2;#3;#4\right)}
\def\emm#1,{{\em #1}}
\def\ba#1,{\overline{#1}}

\def\gen#1,{\langle #1 \rangle}
\def\geno#1,{\langle #1 \rangle_{\infty}}

\def\newop#1{\expandafter\def\csname #1\endcsname{\mathop{\rm
#1}\nolimits}}

\newop{span}
\newop{pr}
\newop{min}
\newop{max}
\newop{deg}
\newop{des}
\newop{Con}
\newop{Ann}
\newop{dim}
\newop{Int}
\newop{Im}
\newop{ker}
\newop{mult}
\newop{rank}
\newop{sign}
\newop{sgn}
\newop{TrCl}
\newop{inf}
\def\cle{\langle\epsilon\rangle}

\begin{document}
\title[Finite automata and pattern avoidance]{
Finite automata and pattern avoidance in words}
\maketitle

\begin{center}
Petter Br\"and\'en\\
Matematik, Chalmers tekniska h\"ogskola och G\"oteborgs
universitet\\
S-412~96  G\"oteborg, Sweden\\[-0.4ex]
{\it branden@math.chalmers.se}\\[1.8ex]
Toufik Mansour\\
Department of Mathematics, Haifa University\\
31905 Haifa, Israel\\[-0.4ex]
 {\it toufik@math.haifa.ac.il}
\end{center}
%=================================================================
\begin{abstract}
We say that a word $w$ on a totally ordered alphabet avoids the word
$v$ if there are no subsequences in $w$ order-equivalent to $v$.
In this paper we suggest a new approach to the enumeration of
words on at most $k$ letters avoiding a given pattern.
By studying an automaton which for fixed $k$ generates the words
avoiding a given pattern
we derive several previously known results for these kind of
problems, as well as many new. In particular,
we give a simple proof of the formula \cite{Reg1998} for
exact asymptotics for the number of words on
$k$ letters of length $n$ that avoids the pattern  $12\cdots(\ell+1)$.
Moreover, we give the first combinatorial proof of the exact
formula \cite{Burstein} for the number of words on $k$ letters of length
$n$ avoiding a three letter permutation pattern.
\end{abstract}
\noindent {\sc 2000 Mathematics Subject Classification}: 05A05,
05A15, 68Q45.
%=================================================================
\section{Introduction}
In this paper we study pattern avoidance in words. The subject of pattern
avoidance in permutations has thrived in the last decades, see \cite{wilf} and
the references there. Only very recently Alon and Friedgut \cite{AloFri2000}
studied pattern avoidance in words to achieve an upper bound on the number
of permutations in $S_n$ avoiding a given pattern. We study pattern
avoidance in words by defining a finite
automaton that generates the words avoiding a given pattern and use the
transfer matrix method to count them. By this approach we are
able to find the asymptotics, as $n \rightarrow \infty$, for the number of
words on $k$ letters of length $n$ avoiding a pattern $p$, as well
as exact enumeration results. In particular we re-derive
Regev's \cite{Reg1998}
result on the exact asymptotics for the number of
words on $k$ letters of length $n$ avoiding a pattern $12 \cdots (\ell+1)$,
and give the first combinatorial proof of a formula
for the number of
words on $k$ letters of length $n$ avoiding the pattern $123$.

Let $S_n$ denote the set of permutations of the set
$[n]:=\{1,2,\ldots,n\}$. If $\sigma \in S_k$ and $\tau \in S_n$,
we say that $\tau$ {\em contains} $\sigma$ if there is a sequence
$1 \leq t_1 < t_2 < \cdots < t_k \leq n$ of integers such that for
all $1 \leq i,j \leq k$ we have $\tau(t_i)  \leq \tau(t_j)$ if and
only if $\sigma(i) \leq \sigma(j)$. Here $\sigma$ is called a {\em
pattern}. If $\tau$ does not contain $\sigma$ we say that $\tau$
{\em avoids} $\sigma$. In the study of pattern avoidance the focus
has been on enumerating and giving estimates to the number of
elements in the set $S_n(\sigma)$, the set of permutations in
$S_n$ that avoids $\sigma$. Maybe the most interesting open
problem in the field is: Does there exists a constant $c$ such
that $|S_n(\tau)|<c^n$ for all $n \geq 0$? This problem is
equivalent to the seemingly stronger statement, see
\cite{Arratia}:
\begin{conjecture}\label{conSW} {\rm ({\bf Stanley, Wilf})}
For any pattern $\tau\in S_\ell$, the limit
$$\lim_{n\rightarrow\infty} |S_n(\tau)|^{\frac{1}{n}},$$
exists and is finite.
\end{conjecture}
The conjecture has been verified for {\em layered} patterns \cite{Bon1999}
and for all patterns  which
can be written as an increasing subsequence followed by a
decreasing \cite{AloFri2000}. In the latter reference Alon and Friedgut
proved a weaker version of Conjecture \ref{conSW}, namely:
For any permutation $\sigma$ there
exists a constant $c=c(\sigma)$ such that $|S_n(\sigma)|\leq
c^{n\gamma^{\star}(n)}$, where $\gamma^{\star}$ is an extremely
slow growing function, related to the Ackermann hierarchy. The method of
proof in \cite{AloFri2000} was by considering pattern avoidance in words.
This is also the theme of this paper.

Denote by $[k]^*$ the set of all finite words with letters in
$[k]$.
If $w=w_1w_2 \cdots w_s \in [k]^*$ and $v = v_1v_2 \cdots
v_r \in [m]^*$ where $r \leq s$, we say that $w$ {\em contains} the
{\em pattern} $v$ if
there is a sequence $1 \leq t_1<t_2< \cdots <t_r \leq s$ such that
for all $1 \leq i,j \leq s$ we have
$$
w_{t_i} \leq w_{t_j} \ \  \mbox{ if and only if } \ \ v_i \leq
v_j.
$$
If $w$ does not contain $v$ we say that $w$ {\em avoids} $v$.
For example, the word $w=323122411\in [4]^9$ avoids the pattern $132$
and contains the patterns $123$, $212$, $213$, $231$, $312$, and $321$. If
$S$ is any set of finite words we denote the set of
words in $S$ that avoids $v$ by $S(v)$.

The history of pattern avoidance in words is not as rich as the one
in permutations. We mention the references
\cite{Atkinson,AloFri2000,Burstein,BM,Klazar,Reg1998}.
In \cite{Reg1998} Regev gave a complete answer for the asymptotics
for $|[k]^n(p_\ell)|$ when $n \rightarrow \infty$, where
$p_\ell=12\cdots(\ell+1)$:
% Besides \cite{AloFri2000},
% \mps{Klazar}
% up to now, forbidden patterns in words dealt with the cases where
% $p$ is a permutation of length $3$ (see Burstein \cite{Burstein}
% and Atkinson et al. \cite{Atkinson}), and $p$ contain repeated
% letters of length $3$ (see Burstein and Mansour \cite{BM}). In the
% case $p_\ell=12\ldots(\ell+1)$, Regev~\cite{Reg1998} gave a complete
% answer for the asymptotics of $[k]^n(p_\ell)$ when $n \rightarrow \infty$.
\begin{theorem}[Regev]\label{regevthm}
For all $k \geq \ell$ we have
$$|[k]^n(p_\ell)| \simeq C_{\ell,k}n^{\ell(k-\ell)}\ell^n \ \ \ \ (n
\rightarrow \infty),$$ where
$$C_{\ell,k}^{-1}= \ell^{\ell(k-\ell)}
\displaystyle\prod_{i=1}^\ell\prod_{j=1}^{k-\ell}(i+j-1).$$
\end{theorem}
\subsection{Organization of the paper}
The paper is organized as follows. In Section \ref{def} we present
the relevant definitions and attain some preliminary results, and
in Section \ref{transfer} we use the transfer matrix method to
determine the asymptotic growth for the sequence $n \mapsto
|[k]^n(p)|$. In Section \ref{inc} we study the special features of
the automaton, $\Au(p_\ell,k)$, which generates the words with
letters in $[k]$ that avoids the increasing pattern $12\cdots
(\ell+1)$. Here we will give a simple proof of Theorem
\ref{regevthm} using the transfer matrix method and give a
combinatorial proof for the formula \cite{Burstein} for
$|[k]^n(p)|$, where $p$ is any permutation pattern of length
three. We also consider the diagonal sequence $|[n]^n(123)|$ and
determine its asymptotic growth as well as showing that its
generating function is transcendental. We conclude the paper by
indicating further problems connected to the work in this paper.

%=============================
\section{Definitions and preliminary results}\label{def}

Given a word-pattern $p$ and an integer $k>0$ we define an equivalence
relation $\sim_p$ on
$[k]^*$ by: $v \sim_p w$ if for all words $r \in [k]^*$ we have
$$vr\mbox{ avoids } p   \ \ \mbox{ if and only if  } \ \ wr \mbox{ avoids } p.$$
For example, if $p=132$, $k\geq4$, $v=13$ and $w=14$, then $v \nsim_p
w$, since $133$ avoids $p$ but $143$ contains $p$.
At first sight it may seem difficult to determine if $v \sim_p w$, since
a priori there is an infinite number of right factors $r$ to check.
By the following
lemma we have to check only a finite number words $r$.

\begin{lemma}\label{pro1}
Let $p$ be a pattern of length $\ell$ and let $v,w\in [k]^*$ be any
two words. Then $v \sim_p w$ if and only if for all words $r \in [k]^s$,
$0 \leq s \leq \ell$, we have
$$
vr \mbox{ avoids } p   \ \ \mbox{ if and only if  } \ \ wr \mbox{
avoids } p.
$$
\end{lemma}
\begin{proof}
Define an equivalence relation $\sim_p'$ on $[k]^*$ by:
$v \sim_p' w$ if for all words $r \in [k]^s$,
$0 \leq s \leq \ell$, we have
$$
vr \mbox{ avoids } p   \ \ \mbox{ if and only if  } \ \ wr \mbox{
avoids } p.
$$
Clearly, $v \sim_p w$ implies $v \sim_p' w$. On the other hand if
$v \nsim_p w$ we may assume that there is an $r \in [k]^*$ such that
$vr$ contains $p$ and $wr$ avoids $p$. Any occurrence of $p$ in $vr$ can use
at most $\ell$ letters of $r$. Thus there is a subsequence $r'$ of
$r$ of length at most $\ell$ such that $vr'$ contains $p$ and $wr'$
avoids $p$, i.e., $v \nsim_p' w$.
\end{proof}

Let $\E(p,k)$ be the set of equivalence classes of $\sim_p$. By
Lemma \ref{pro1} the number of equivalence classes is finite. We
denote the equivalence class of a word $w$ by $\langle w \rangle$.
The equivalence classes of $\sim_p$ for $p\in S_3$
and $k=3,4,5$ are given in Table~1.

\begin{definition}\label{defauto}
Given a positive integer $k$ and a pattern $p$ we define a {\em
finite automaton}\footnote{For a definition of a finite automaton,
see \cite{ASU} and references therein.}, $\Au(p,k)=(\E(p,k), [k],
\delta,\langle \e \rangle, \E(p,k) \setminus \{\langle p
\rangle\})$, by
\begin{itemize}
\item the {\em states} are, $\E(p,k)$, the equivalence-classes of $\sim_p$,
\item $[k]$ is the {\em input alphabet},
\item $\delta : \E(p,k) \times [k] \rightarrow \E(p,k)$ is the
{\em transition function} defined by
$\delta(\langle w \rangle,i)=\langle wi \rangle$, where $wi$ is $w$
concatenated with the letter $i \in [k]$,
\item $\langle \e \rangle$ is the {\em initial state},
      where $\e$ is the empty word,
\item all states but $\langle p \rangle$ are {\em final states}.
\end{itemize}
\end{definition}
For an example see Fig. \ref{fig2314}.
\begin{center}
\begin{table}\label{tabequi}
{\footnotesize
\begin{tabular}{|c|c|l|}\hline
   $k$ & $p$   & The equivalences classes in $\E(p,k)$ \\ \hline
   $3$ & $123$ & $\cle, \langle1\rangle, \langle12\rangle, \langle123\rangle$  \\
       & $132$ & $\cle, \langle1\rangle, \langle13\rangle, \langle132\rangle$  \\
       & $213$ & $\cle, \langle2\rangle, \langle21\rangle, \langle213\rangle$  \\
       & $231$ & $\cle, \langle2\rangle, \langle23\rangle, \langle231\rangle$  \\
       & $312$ & $\cle, \langle3\rangle, \langle31\rangle, \langle312\rangle$  \\
       & $321$ & $\cle, \langle3\rangle, \langle32\rangle, \langle321\rangle$
       \\\hline\hline

   $4$ & $123$ & $\cle, \langle1\rangle, \langle2\rangle, \langle12\rangle, \langle13\rangle, \langle23\rangle, \langle123\rangle$ \\
       & $132$ & $\cle, \langle1\rangle, \langle2\rangle, \langle13\rangle, \langle14\rangle, \langle24\rangle, \langle132\rangle, \langle241\rangle$ \\
       & $213$ & $\cle, \langle2\rangle, \langle3\rangle, \langle21\rangle, \langle23\rangle, \langle31\rangle, \langle32\rangle, \langle213\rangle$ \\
       & $231$ & $\cle, \langle2\rangle, \langle3\rangle, \langle23\rangle, \langle24\rangle, \langle32\rangle, \langle34\rangle, \langle231\rangle$ \\
       & $312$ & $\cle, \langle3\rangle, \langle4\rangle, \langle31\rangle, \langle41\rangle, \langle42\rangle, \langle312\rangle, \langle314\rangle$ \\
       & $321$ & $\cle, \langle3\rangle, \langle4\rangle, \langle32\rangle, \langle42\rangle, \langle43\rangle, \langle321\rangle$
       \\\hline\hline

   $5$ & $123$ & $\cle, \langle1\rangle, \langle2\rangle, \langle3\rangle, \langle12\rangle, \langle13\rangle, \langle14\rangle, \langle23\rangle, \langle24\rangle, \langle34\rangle, \langle123\rangle$ \\
       & $132$ & $\cle, \langle1\rangle, \langle2\rangle, \langle3\rangle, \langle13\rangle, \langle14\rangle, \langle15\rangle, \langle24\rangle, \langle25\rangle, \langle35\rangle, \langle132\rangle, \langle241\rangle, \langle251\rangle, \langle351\rangle, \langle352\rangle, \langle3513\rangle$ \\
       & $213$ & $\cle, \langle2\rangle, \langle3\rangle, \langle4\rangle, \langle21\rangle, \langle23\rangle, \langle24\rangle, \langle31\rangle, \langle32\rangle, \langle34\rangle, \langle41\rangle, \langle42\rangle, \langle43\rangle, \langle213\rangle, \langle234\rangle, \langle243\rangle$ \\
       & $231$ & $\cle, \langle2\rangle, \langle3\rangle, \langle4\rangle, \langle23\rangle, \langle24\rangle, \langle25\rangle, \langle32\rangle, \langle34\rangle, \langle35\rangle, \langle42\rangle, \langle43\rangle, \langle45\rangle, \langle231\rangle, \langle243\rangle, \langle432\rangle$ \\
       & $312$ & $\cle, \langle3\rangle, \langle4\rangle, \langle5\rangle, \langle31\rangle, \langle41\rangle, \langle42\rangle, \langle51\rangle, \langle52\rangle, \langle53\rangle, \langle312\rangle, \langle314\rangle, \langle315\rangle, \langle415\rangle, \langle425\rangle, \langle3153\rangle$ \\
       & $321$ & $\cle, \langle3\rangle, \langle4\rangle, \langle5\rangle, \langle32\rangle, \langle42\rangle,
       \langle43\rangle, \langle52\rangle, \langle53\rangle, \langle54\rangle, \langle321\rangle$\\ \hline
\end{tabular}}
\caption{The equivalence classes of $\sim_p$ for $p\in S_3$ and
$k=3,4,5$.}
\end{table}
\end{center}

We will identify $\Au(p,k)$ with the (labelled) directed graph
with vertices $\E(p,k)$ and with a (labelled) edge
$\stackrel{i}{\longrightarrow}$ between $\langle v \rangle$ and
$\langle w \rangle$ if $vi \sim_p w$. Clearly, we may order the
states as $x_1,x_2, \ldots, x_e$ so that if $i<j$ there is no path
from $x_j$ to $x_i$. The {\em transition matrix}, $T(p,k)$, of
$\Au(p,k)$ is the matrix of size $e \times e$ with non-negative
integer coefficients defined by:
        $$[T(p,k)]_{ij}= |\{s \in [k] : \delta(x_i,s)=x_j\}|.$$
Thus $[T(p,k)]_{ij}$ counts the number of edges between $x_i$ and
$x_j$, and $T(p,k)$ is triangular.

\begin{example}\label{ex2314}
If $p=2314$ and $k=5$, then it is easy to check (see~\cite{M})
that the states are $\cle$, $\langle2\rangle$, $\langle3\rangle$,
$\langle32\rangle$, $\langle34\rangle$, $\langle24\rangle$,
$\langle23\rangle$, $\langle324\rangle$, $\langle341\rangle$,
$\langle241\rangle$, $\langle234\rangle$, $\langle2342\rangle$,
$\langle231\rangle$, and $\langle2314\rangle$ (see Fig. \ref{fig2314}).
\begin{figure}[h]
%\hspace*{-100pt}
\epsfxsize=400.0pt \epsffile{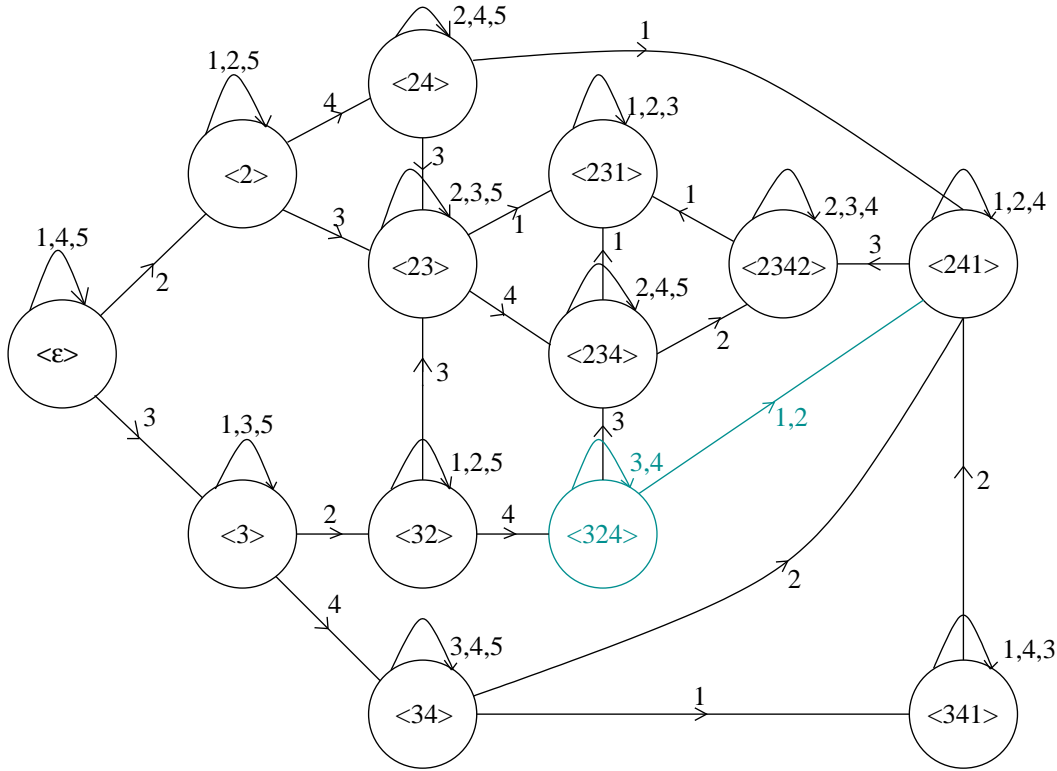}
\caption{The figure shows the final states in the automaton $\Au(2314,5)$.}
\label{fig2314}
\end{figure}

Note that there are two edges between the states
$\langle324\rangle$ and $\langle241\rangle$, namely
$\langle324\rangle\stackrel{1}{\longrightarrow}\langle241\rangle$
and
$\langle324\rangle\stackrel{2}{\longrightarrow}\langle241\rangle$.
Moreover, all final states in $\Au(2314,5)$
have $3$ loops, except $\langle 324\rangle$ which has $2$ loops.
\end{example}

The following simple lemma will be helpful in finding the asymptotic growth
of the sequence $|[n]^k(p)|$, for fixed $k$.

\begin{lemma}\label{loops}
Let the automaton $\Au(p,k)$ be given, let $d$ be the number of
distinct letters in $p$ and suppose that $k \geq d-1$.
If $\langle v \rangle$ is any state different from
$\langle p \rangle$, then the number of loops at $\langle v
\rangle$ does not exceed $d-1$. Moreover,
there are exactly $d-1$ loops at $\langle \e \rangle$.
\end{lemma}
\begin{proof}
Suppose that there are more than $d-1$ loops at $\langle v \rangle$. Then
the loops use at least $d$ different labels. From
these labels we can form a word $w$ order-isomorphic to $p$. But
then $vw \sim_p v$ which is a contradiction.

Let $p_1$ be the first letter of $p$. Then, if
$i < p_1$ or $i > k-d+p_1$ we have $i \sim_p \e$. But there are
$d-1$ such $i$s, which proves the lemma.
\end{proof}

Although pattern avoidance in words and pattern avoidance in permutations
share many common features, there are some important aspects in which they
differ.
For permutations there are three simple operations, $f$, that respects
pattern-avoidance in the sense that $f(\tau)$ avoids $f(\sigma)$ if and
only if $\tau$ avoids $\sigma$, namely the reversal, the complement
and the inverse of a permutation. The first two operations have
obvious generalizations to words, while the inverse does not. It
has in fact been an open question to construct an inverse for words
possessing ``the right'' properties. Such an inverse was recently
constructed by Hohlweg and Reutenauer \cite{HohlwegReutenauer}.
Unfortunately it is not possible to construct an inverse
that respects pattern avoidance in words, which would imply the identity
$
|[k]^n(p)| = |[k]^n(p^{-1})|,
$
for all $k,n \geq 0$ and permutation patterns $p$.
The first counter example to this is
$|[5]^7(1342)|=67854 > 67853=|[5]^7(1423)|$, see Table \ref{tab4}.
If $w \in [k]^n$ let the complement of $w$ in $[k]^n$ be
$w^c= (k+1-w_1)(k+1-w_2)\cdots (k+1-w_n)$. Then we have in fact that
$\Au(p,k)$ and $\Au(p^c,k)$ are isomorphic as automata for any $p \in [k]^*$,
since $v \sim_p w$ if and only if $v^c \sim_{p^c} w^c$.
% \begin{lemma}\label{compl}
% Let $p=p_1p_2\ldots p_d$ be any pattern and $k\geq1$ with
% automaton $$\Au(p,k)=(\E(p,k), [k], \delta,\langle \e \rangle,
% \E(p,k) \setminus \{\langle p \rangle\}).$$ Then, the automaton
% $\Au(p',k)$, where $p'=(k+1-p_1)(k+1-p_2)\ldots(k+1-p_d)$, is
% given by
% \begin{itemize}
% \item the vertices are $\langle(k+1-w_1)\ldots(k+1-w_s)\rangle$, where
% $\langle w_1\ldots w_s\rangle\in\E(p,k)$,

% \item the edges are $\stackrel{k+1-i}{\longrightarrow}$ between two
% vertices $\langle (k+1-w_1)\ldots (k+1-w_s)\rangle$ and $\langle
% (k+1-u_1)\ldots (k+1-u_t)\rangle$, where
% $\stackrel{i}{\longrightarrow}$ is an edge in $\Au(p,k)$ between
% the vertices $\langle w_1w_2\ldots w_s\rangle$ and $\langle
% u_1u_2\ldots u_t\rangle$.
% \end{itemize}
% \end{lemma}
% \begin{proof}
% Let us denote the equivalence relation $\sim_p$ on $[k]^*$ for given
% a pattern $p$ by $\sim_p$. Hence, by definitions it is easy to see
% that $\E(p,k)$ the equivalence classes of $\sim_p$ if and only if
% $\E(p',k)$ the equivalence classes of $\sim_{p'}$. Moreover, there
% is an edge $\stackrel{i}{\longrightarrow}$ between two vertices
% $\langle w_1w_2\ldots w_s\rangle$ and $\langle u_1u_2\ldots
% u_t\rangle$ in $\Au(p,k)$ if $v_1\ldots v_s i\sim_p w_1\ldots w_t$
% if and only if
% $(k+1-v_1)\ldots(k+1-v_s)(k+1-i)\sim_{p'}(k+1-w_1)\ldots(k+1-w_t)$,
% equivalently $\stackrel{k+1-i}{\longrightarrow}$ is an edge in
% $\Au(p',k)$ between the vertices $\langle (k+1-w_1)\ldots
% (k+1-w_s)\rangle$ and $\langle (k+1-u_1)\ldots (k+1-u_t)\rangle$.
% \end{proof}
Certainly $w$ avoids $p$ if and only if $w^r$ avoids
$p^r$, where $^r$ is the reversal operator and $w$ and $p$ are any words.
However $\Au(p,k)$ and $\Au(p^r,k)$ are not in general isomorphic. Indeed,
for $p=2314$ and $k=5$ we have that $|\E(2314,5)|=13$ and $|\E(4132,5)|=14$.
%===========================================================
\section{Transfer matrix method}\label{transfer}
In this section we use the transfer matrix method (see
\cite[Theorem~4.7.2]{Stanley1})
to obtain information about the sequences $|[k]^n(p)|$.
Given a matrix $A$ let $(A;i,j)$ be the matrix with row $i$ and
column $j$ deleted. If $p$ is a pattern and $k>0$ an  positive integer let
$T'(p,k)= (T(p,k),e_k-1,e_k-1)$.

\begin{theorem}\label{gf} Let $k$
be a positive integer, $p$ be a pattern and $e_k$ be the
number of states in $\Au(p,k)$. Then the generating function for
$|[k]^n(p)|$ is
$$
\sum_{n \geq 0} |[k]^n(p)|x^n=
\frac {\sum_{j=1}^{e_k-1}{(-1)^{j+1}\det(I-xT',j,1)}}
{\prod_{i=1}^{e_k-1}(1-\lambda_ix)}= \frac{ \det B(x)}{\prod_{i=1}^{e_k-1}(1-\lambda_ix)},
$$
where $\lambda_i$ is the number of loops at state $x_i$, and $B(x)$ is the
matrix obtained by replacing the first column in $I-xT'$ with a column
of all ones.
\end{theorem}
\begin{proof}
The theorem follows from the transfer matrix method, see
\cite[Theorem 4.7.2]{Stanley1}, since we want to count the number
of paths from $\langle \e \rangle$ to any state other than
$\langle p \rangle$ of length $n$ in $\Au(p,k)$.
\end{proof}

Regev \cite{Reg1998} computed the exact asymptotics for
$|[k]^n(p_\ell)|$, where $p_\ell$ is the increasing pattern $12
\cdots(\ell+1)$ and $n \to \infty$. We will next
find the exact asymptotics (up to a constant) for $|[k]^n(p)|$ for
all patterns $p$. Given two sequences $\{a_n\}$ and $\{b_n\}$ of
real numbers, we denote $a_n \simeq b_n$ if $\lim_{n
\rightarrow\infty}\frac{a_n}{b_n} = 1$. A path in $\Au(p,k)$ is
called {\em simple} if it starts at $\langle \e \rangle$, does not use any
loops, and does not end in $\langle p \rangle$.

\begin{theorem}\label{asymptotics}
Let $p$ be any pattern with $d$ distinct letters and let $k \geq
d-1$ be given. Then there is a constant $C>0$ such that
$$
|[k]^n(p)| \simeq Cn^M(d-1)^n \ \ \ \ (n \rightarrow \infty),
$$
where $M+1$ is the maximum number of states with $d-1$ loops, in
a simple path.
\end{theorem}
\begin{proof}
Let $P :=x_1, x_2, \ldots, x_j$ be a simple path in $\Au(p,k)$.
Moreover, let $\ell_j$ be the number of loops at state
$x_j$. Then $|[k]^n(p)|=\sum_P N(P,n)$ where
$$
N(P,n)= \sum_{\alpha_1+\cdots+\alpha_j=n-j+1}\ell_1^{\alpha_1}
\ell_2^{\alpha_2} \cdots \ell_j^{\alpha_j},
$$
and the sum is over all weak compositions of $n-j+1$ into at most
$j$ parts. Now, $N(P,n)$ is equal to the coefficient to
$t^{n-j+1}$ in $(1-\ell_1t)^{-1}\cdots(1-\ell_jt)^{-1}$. Let $r$ be
the number of $i$ such that $\ell_i = d-1$. Note that by Lemma \ref{loops}
$r$ is greater than or equal to one. The dominant term
of $(1-\ell_1t)^{-1}\cdots(1-\ell_jt)^{-1}$ is (by partial fraction
decomposition) equal to
            $$\frac {f(t)}{(1-(d-1)t)^r},$$
where $f(t)$ is a polynomial of degree less than $r$ and
$f((d-1)^{-1}) \neq 0$. By well known results it follows that
$N(P,n) \simeq C(P)(d-1)^n n^{r-1}$, where $C(P) > 0$ is a
constant depending on $P$ and $k$. Taking the greatest possible
$r$ yields the desired results.
\end{proof}

When there are exactly $d-1$ loops at every state except $\langle
p \rangle$ in $\Au(p,k)$, then it follows from Theorem \ref{gf}
that $|[k]^n(p)|= (d-1)^n Q(n)$, where $Q$ is a polynomial in $n$.
We have in fact:

\begin{corollary}\label{pollypolly}
Let $\Au(p,k)$ be such that all states but $\langle p \rangle$ has
exactly $d-1$ loops. Then
$$
|[k]^n(p)|= \sum_{j=0}^M a_j(d-1)^{n-j} \binom n j,
$$
where $a_j$ counts the number of simple paths of length $j$ in $\Au(p,k)$.
Moreover, if $p$ is a pattern of length
$\ell+1$ then $a_j = (k-d+1)^j$ for all $j=0,1,\ldots,\ell$.
\end{corollary}
\begin{proof}
The corollary follows from the proof of Theorem \ref{asymptotics}
since $N(P,n)=(d-1)^{n-j}\binom n j$. If $p$ is a pattern of
length $\ell + 1$ then we have that $a_j = (k-d+1)^j$ where
$j=0,1,\ldots,\ell$,  since $k^j=\sum_{i=0}^ja_i(d-1)^{j-i}\binom
j i$ for all $j=0,1,\ldots,\ell$.
\end{proof}
As an example of Corollary \ref{pollypolly} we note that
if $p$ is any pattern of length
$\ell+1$ with exactly $d$ different letters then
$$
|[d]^n(p)|=\sum_{j=0}^{\ell}(d-1)^{n-j}\binom{n}{j}.
$$
% \subsection{Pattern starting with $1$}
% If $p$ is any pattern let $1p$ be the pattern obtained by adding
% one to every letter of $p$ and concatenating a one to the left.

% \begin{lemma}\label{addone}
% Suppose that $p$ is a pattern with first entry equal to one,
% satisfying
% $$
% |[k]^n(p)| \leq C^kD^n ,
% $$
% for some constants $C,D >0$ and all $n,k \geq 0$. Then
% $$
% |[k]^n(1p)| \leq (2C)^k(2D+2)^n.
% $$
% \end{lemma}
% \begin{proof}
% Let $w \in [k]^n(1p)$. The letters of $w$ which are the first
% letter of an occurrence of $p$ must necessarily be weakly
% decreasing. The remaining letters avoids $p$. Thus
% \begin{eqnarray*}
% |[k]^n(1p)| &\leq& \sum_{i=0}^n \binom n i \binom {i+k-1}{k-1} C^kD^{n-i} \\
%             &\leq& 2^{n+k-1}C^k(D+1)^n,
% \end{eqnarray*}
% since there are $\binom {i+k-1}{k-1}$ weakly decreasing sequences
% of length $i$ using at most $k$ letters. This proves the lemma.
% \end{proof}

% Let $p_{\ell} := 12 \cdots \ell (\ell+1)$ be the increasing
% pattern of length $\ell +1$.

% \begin{theorem}\label{thinc1}
% For all $n\geq0$,
%     $$|[k]^n(p_\ell)| \leq  2^{\ell k}(2\ell)^n.$$
% \end{theorem}
% \begin{proof}
% Follows by induction from Lemma \ref{addone}.
% \end{proof}
% We remark that, Theorem~\ref{thinc1} gives that
% $$|S_n(p_\ell)|\leq |[n]^n(p_\ell)|\leq
% (\ell\cdot 2^{\ell+1})^n.$$

\section{The increasing patterns}
We will in this section investigate the properties of $\Au(p_{\ell},k)$,
where $p_\ell = 12\cdots (\ell+1)$. The following lemma describes the
structure of $\Au(p_{\ell},k)$:
\begin{lemma}\label{inc}
Let $k \geq \ell$ be given.
For any subset $S$ of $[k]$ of size $\ell$ let $w_S$ be
the word consisting of the elements of $S$ listed in increasing
order. Then the words $w_S$ together with $p_\ell$
constitute a complete set of representatives for the
equivalence-classes $\E(p_\ell,k)$. In particular we have:
$$
|\E(p_\ell,k)| = \binom k {\ell} + 1.
$$
If $S=\{s_1< \cdots < s_\ell\} \subseteq [k]$ and $j \in [k]$ let
$$S^j = \{s_1< \cdots <s_{i-1}<j<s_{i+1}<\cdots<s_\ell\},$$ where $i$ is the
integer such that $s_{i-1}<j\leq s_i$, ($s_0:=0, s_{\ell+1}:=k+1$). Then
$$
\delta(\langle w_S \rangle,j)=
\begin{cases}
\langle w_{S^j} \rangle \mbox{ if } j \leq s_\ell, \\
\langle p_\ell \rangle \mbox{ otherwise }.
\end{cases}
$$
In particular, the loops of $w_S$ are the elements of $S$.
\end{lemma}
\begin{proof}
It is clear that the words $w_S$ are representatives for different
classes. Let $v \in [k]^*(p_\ell)$. We say that an increasing subword
$x_1x_2\cdots x_j$ of $v$ is {\em extendible} if $x_j \leq  k+j-\ell-1$,
i.e., if we may extend $x_1x_2\cdots x_j$ to an occurrence of $p_\ell$
using letters from $[k]$. Suppose that the maximum length of an
extendible increasing subsequence in $v$ is equal to $s$, $s\leq \ell$.
For $1 \leq j \leq s$ let
$$
r_j(v) := \min\{ x_j : x_1x_2\cdots x_j
\mbox{ is an extendible subword of } v \}.
$$
Clearly $r_1(v)<r_2(v)< \cdots < r_s(v)$. Let
$$
S = \{ r_1(v), r_2(v), \ldots, r_s(v), k+s+1-\ell,k+s-\ell,\ldots,
k-1 ,k\}.
$$
Then we see that $w_S \sim v$.
The statement about the transition function follows from the construction.
\end{proof}
In the sequel we will use some standard notation from the theory
of partitions and symmetric functions. For undefined terminology we refer
the reader to Chapter 7 of \cite{Stanley2}.
\begin{theorem}\label{orderideals}
Define a partial order on the final states in $\Au(p_\ell,k)$ by:
$x \leq y$ if there exists a path from $x$ to $y$ in  $\Au(p_\ell,k)$.
Then this partial order is isomorphic to
$$
J([\ell]\times [k-\ell]),
$$
the lattice of order ideals of the poset $[\ell]\times [k-\ell]$.
\end{theorem}
\begin{proof}
Let $S=\{s_1<s_2< \cdots <s_\ell\}$ and $T=\{t_1<t_2< \cdots <t_\ell\}$ be
subsets of $[k]$.
We claim that there exists a path from $\langle w_S \rangle$ to
$\langle w_T \rangle$ if and only if $s_i \geq t_i$ for all
$1 \leq i \leq \ell$. From this the theorem follows since the latter
poset is isomorphic to the interval $[\emptyset,\lambda_{\ell,k-\ell}]$,
in the Young's lattice, where
$\lambda_{\ell,k-\ell}:= (k-\ell,k-\ell,\ldots, k-\ell)$ is of length $\ell$.
Indeed, consider the bijection defined by:
$$
(s_1,s_2, \ldots, s_\ell) \mapsto
(s_\ell -\ell, s_{\ell-1}- \ell+1, \ldots , s_1-1) \in
[\emptyset,\lambda_{\ell,k-\ell}].
$$
Then $s_i \geq t_i$ for all $1 \leq i \leq j$ if and only if
the image of $S$ is greater than the image of $T$ in
$[\emptyset,\lambda_{\ell,k-\ell}]$. But $[\emptyset,\lambda_{\ell,k-\ell}]$ is
its own dual, so the statement follows from the simple fact that
$[\emptyset,\lambda_{\ell,k-\ell}]$ is
isomorphic to $J([\ell]\times [k-\ell])$.

If there is an edge between  $\langle w_S \rangle$
and $\langle w_T \rangle$, we are done by Lemma \ref{inc}. The ``only if''
direction thus follows by induction on the length of the path.

Now, if $s_i \geq t_i$ for all $1 \leq i \leq \ell$ consider the path
$$
\langle w_S \rangle \stackrel{t_1}{\longrightarrow}
\langle w_{S}t_1 \rangle \stackrel{t_2}{\longrightarrow}
\langle w_{S}t_1t_2 \rangle \stackrel{t_3}{\longrightarrow} \cdots
\stackrel{t_\ell}{\longrightarrow}  \langle w_{S}t_1t_2\cdots t_\ell \rangle .
$$
It is not hard to see that
$\langle w_{S}t_1t_2\cdots t_\ell \rangle = \langle w_T \rangle$, which
completes the proof.
\end{proof}
We now have a different proof of the following theorem of Regev
\cite{Reg1998}:
\begin{theorem}[Regev]
For all $k \geq \ell$ we have
$$
|[k]^n(p_\ell)| \simeq C_{\ell,k}n^{\ell(k-\ell)}\ell^n
\ \ \ \ (n \rightarrow \infty),
$$
where
$$
C_{\ell,k}^{-1}= \ell^{\ell(k-\ell)}
\displaystyle\prod_{i=1}^\ell\prod_{j=1}^{k-\ell}(i+j-1).
$$
\end{theorem}
\begin{proof}
By Corollary \ref{pollypolly} and Theorem \ref{orderideals} we have that
$$
|[k]^n(p_\ell)| \simeq a_M\ell^{-M}\binom n M \ell^n \simeq
                \frac{a_M}{M!}\ell^{-M}n^M \ell^n \ \ \ \
(n \rightarrow \infty),
$$
where $M=\ell(k-\ell)$ and $a_M$ is equal to
the number
of maximal chains in $J([\ell]\times [k-\ell])$. By
\cite[Proposition 7.10.3]{Stanley2} and the hook-length formula
\cite[Corollary 7.21.6]{Stanley2} we have that
$$
a_{\ell(k-\ell)} = f^{\lambda_{\ell,k-\ell}}=\frac {(\ell(k-\ell))!}{
\displaystyle\prod_{i=1}^\ell\prod_{j=1}^{k-\ell}(i+j-1)},
$$
from which the theorem follows.
\end{proof}
It should be clear from the correspondence in Theorem
\ref{orderideals} that the simple paths of length $r$ in
$\Au(p_\ell,k+\ell)$ are in a one-to-one correspondence with
tableaux $T$ of the following type:
\begin{itemize}
\item[(i)] $T$ is weakly increasing in rows and columns,
\item[(ii)] no integer appears in more than one row,
\item[(iii)] the entries of $T$ are exactly $[r]$,
\item[(iv)] the shape of $T$ is confided in $\lambda_{\ell,k}$.
\end{itemize}
Recall that the tableaux satisfying (i) and (ii) above are the
{\em border-strip} tableaux (or  {\em rim-hook} tableaux) of
height zero. We call these tableaux {\em segmented}. Let $a(\ell,k,r)$
denote the number of segmented tableaux satisfying (iii) and (iv),
so that:
\begin{equation}\label{sumbin}
|[k+\ell]^n(p_\ell)| = \sum_{r=0}^{\ell k}\ell^{n-r}a(\ell,k,r)\binom n r.
\end{equation}
The function $a(\ell,k,r)$ is actually a polynomial in $k$ of
degree $r$. To see this let us call a segmented tableau inside
$[\ell] \times [k]$ {\em primitive} if all columns are different,
and let the set of such tableaux of length $i$ with $r$ different
entries be $\PR_{\ell,i,r}$. If we denote the number of elements
in $\PR_{\ell,i,r}$ by $\pr(\ell,i,r)$ we have
$$
a(\ell,k,r) = \sum_{i = r/\ell}^r \pr(\ell,i,r) \binom k i,
$$
since for any such primitive tableaux of length $i$ we may insert
a number $\alpha_1$ copies of the first column before the first
column, a number $\alpha_2$ copies of the second column between
the first and the second column, and so on. After the last column
we may insert a number $\alpha_{i+1}$ columns of all blanks,
requiring that
$$
\alpha_1 + \alpha_2 + \cdots + \alpha_{i+1} = k-i.
$$
Thus there are $\binom k i$ segmented tableaux arising from a
given primitive one. The numbers $\pr(\ell,i,r)$ are in general
hard to count, but there are two special cases which are nice,
namely $\pr(\ell,r,r)$ and $\pr(2,i,r)$. We start by counting
$\pr(\ell,r,r)$.
\begin{theorem}\label{perms}
With definitions as above:
$$
\pr(\ell,n,n) = |S_n(p_\ell)|.
$$
\end{theorem}
\begin{proof}
We will define a bijection between $S_n$ and
$\cup_{\ell \geq 0}\PR_{\ell,n,n}$ such that the height of the tableau
corresponds to the greatest increasing subsequence in the permutation. Recall
the definition of $r_i(v)$ in the proof of Lemma \ref{inc}, and
let $r(v)=(r_1(v),r_2(v), \ldots, r_\ell(v))$, where $\ell$ is the length of
the longest increasing subsequence in $v$.  Let $k$ be
big enough so that all increasing subsequences in permutations in $S_n$ are
considered extendible.

Now, if $\pi=\pi_1\pi_2 \cdots \pi_n$ is any permutation in $S_n$
define $T=T(\pi)$
as follows. Let the first column of $T$ be $r(\pi)$ the second
column be $r(\pi_1 \cdots \pi_{n-1})$ and so on. The image of the
permutation $351462$ is:
$$
T(351462) =
\begin{array}{cccccccccc}
  1&1&1&1&3&3\\
  2&4&4&5&5&\\
  6&6& & & &
\end{array}.
$$
By
Lemma \ref{inc} we have that $T(\pi) \in \PR_{\ell,n,n}$. Moreover from Lemma
\ref{inc} we also get that a
tableau $T$ is the image of some $\pi \in S_n$ if and only if
\begin{itemize}
\item[(a)] $T$ has $n$ columns and entries $1,2, \ldots, n$,
\item[(b)] Let $T^i$ denote the $i$th column. If $i<j$ then
$T^i$ is smaller than $T^j$ in the product order. (If $T^i$ and $T^j$ have
different size fill the empty slots of $T^j$ with $n+1$),
\item[(c)] Exactly one new entry appears every time you move from
      $T^{i+1}$ to $T^i$.
\end{itemize}
Now, if $T \in \cup_{\ell \geq 0}\PR_{\ell,n,n}$ condition (a) and
(b) are trivially satisfied. At least one new entry appears every
time we move from $T^{i+1}$ to $T^i$, since otherwise $T^i=T^{i+1}$ and
$T$ fails to be primitive. On the other hand if there appear
more than one new entry in a transition then in a later transition
there must appear no new entry, since $T$ has $n$ columns and
$n$ entries. This verifies condition (c) and the theorem follows.
\end{proof}
A special case of Theorem \ref{perms} is that
$\pr(2,n,n)=C_n$, the $n$th Catalan number. This is also
a special case of the next theorem. Note, that Theorem \ref{123-bij}
is what we need to have combinatorial proof of a closed formula, see
Theorem \ref{closed}, for
the numbers $|[k]^n(123)|$. Burstein \cite{Burstein} achieved a different,
but of course equivalent, formula for $|[k]^n(123)|$, but not in a
bijective manner.
\begin{theorem}\label{123-bij}
With definitions as above:
$$
\pr(2,i,r) = \frac 1 {i+1}\binom{2i}i \binom i {r-i}.
$$
\end{theorem}
Before we give a proof of Theorem \ref{123-bij} we will need some
definitions and a lemma. Let $\PR^+(2,s,r)$ be the tableaux in
$\PR(2,s,r)$ that fill up the shape $[2]\times [r]$, and let
$\pr^+(2,s,r):=|\PR^+(2,s,r)|$. Then
$\pr(2,s,r)=\pr^+(2,s,r)+\pr^+(2,s,r+1)$ since we get the tableaux
that do not fill up the shape by deleting all entries $r+1$. To
prove the theorem we will show that $\pr^+(2,s,r)=\binom
{s-1}{2s-r}C_s$, where $C_s$ is the $s$th Catalan number.

We first define an operation $+$ that takes tableaux with $r$
different entries to tableaux with $r+1$ different entries. Let $T
\in \PR^+(2,s,r)$. Suppose that $j$ is an index such that
$T_{ij}=T_{i(j+1)}$ for some $i=1,2$. Write $T$ as $T=LR$ where
$L$ is the $j$ first columns and $R$ is the $s-j$ last columns.
Let $R'$ be the array order equivalent with $R$ with entries
the same as $R$, add $r+1$, take away $T_{i(j+1)}$ (two arrays $A$
and $B$ are said to be order equivalent if $A_{ij} \leq A_{i'j'}$ if and only
if $B_{ij} \leq B_{i'j'}$ for all $i,j,i',j'$). We define $T+j$ to
be the tableaux $T+j:=LR'$. In $T$ there are exactly $t=2s-r$
indices $j \in [s-1]$ such that $T_{ij}=T_{i(j+1)}$ for some
$i=1,2$. Let $S=\{s_1<s_2< \cdots <s_t\}$ be these indices and
define a function $\Phi : \PR^+(2,s,r) \rightarrow \binom {[s-1]}t
\times \ST_{2,s}$, where $\ST_{2,s}$ is the set of standard
tableaux of shape $[2]\times[s]$, by
$$
\Phi(T)=(S,T+s_t+s_{t-1}+\cdots + s_1).
$$
The fact that $\Phi$ is a bijection will prove the theorem, since by
the hook-length formula we have $|\ST_{2,s}|=C_s$.
To find the inverse of $\Phi$ we need a kind of inverse operation to $+$.

Let $T \in \PR^+(2,s,r)$ and $1 \leq b \leq s-1$ be such that
$T_{1b}<T_{1(b+1)}$ and $T_{2b} < T_{2(b+1)}$. Define two arrays
$T|_b$ and $T|^b$ as follows. Write $T = LR$ where $L$ are the $b$ first
columns and $R$ are the $s-b$ last columns.
Define $T|^b:=L'R'$, to be the array where $L=L'$ and $R'$ is the unique
array order equivalent with $R$, with entries the same as $R$ add $T_{1b}$
take away $r$. Similarly, let $T|_b:=L'R'$, be the array with
$L=L'$ and where $R'$ is the unique
array order equivalent with $R$, with entries the same as $R$, add $T_{2b}$
take away $r$.

\begin{eqnarray*}
\left.\begin{array}{cccccc}
 1&2&4&4\\
 3&5&6&7
\end{array}
\right|^2 &=&
\begin{array}{cccccc}
 1&2&2&2\\
 3&5&4&6
\end{array}
\\
\left.\begin{array}{cccccc}
 1&2&4&4\\
 3&5&6&7
\end{array}
\right|_2 &=&
\begin{array}{cccccc}
 1&2&4&4\\
 3&5&5&6
\end{array}
\end{eqnarray*}
Note that exactly one of $T|^2$ and $T|_2$ above is a primitive
segmented tableaux. This is no accident.
\begin{lemma}\label{bijlemma}
Let $T \in \PR^+(2,s,r)$ and $1 \leq b \leq s-1$ be such that
$T_{1b}<T_{1(b+1)}$ and $T_{2b} < T_{2(b+1)}$. Then
\begin{eqnarray*}
T|_b \in \PR^+(2,s,r-1) \ \ &\Leftrightarrow& T|^b \notin \PR^+(2,s,r-1)\\
&\Leftrightarrow& T_{2(b+1)}= T_{2b}+1
\end{eqnarray*}
Moreover, if $B=T|^b \in \PR^+(2,s,r-1)$ then $B_{1b}=B_{1(b+1)}$
and if $A=T|_b \in \PR^+(2,s,r-1)$ then $A_{1b}=A_{1(b+1)}.$
\end{lemma}
\begin{proof}
Consider  $A:=T|_b$.
All entries in $T$ that are smaller than $T_{2b}$ will
be mapped on themselves and $A_{ij}=T_{ij}-1$ for $A_{ij}>T_{2b}$. Therefore
$A \in \PR^+(2,s,r-1)$  if and only if  $T_{2(b+1)}= T_{2b}+1$ (since
otherwise the entry $T_{2b}$ will appear in both the first and
the second row).

Consider $B:=T|^b$. Let $y_i$, $i=1,2, \ldots, h$ be the entries in $T$
satisfying $T_{2b}<y_i\leq T_{2(b+1)}$ ordered by size. Then the entry $y_1$
will be mapped
to an element smaller than $T_{2b}$ and $y_i$ will be
mapped to $y_{i-1}$ for $i>1$. Thus $B \in \PR^+(2,s,r-1)$ if and
only if $T_{2(b+1)}> T_{2b}+1$ as claimed.

The last statement is a direct consequence of the above proof.
\end{proof}
We are now ready to give a proof of Theorem \ref{123-bij}.
\begin{proof}[Proof of Theorem \ref{123-bij}]
If $T \in \PR^+(2,s,r)$ and $1 \leq b \leq s-1$ are such that
$T_{1b}<T_{1(b+1)}$ and $T_{2b} < T_{2(b+1)}$ we
define $T- b$ to be the one of the arrays $T|_b$ and $T|^b$ which
is in $\PR^+(2,s,r-1)$. By Lemma \ref{bijlemma} we have that
\begin{eqnarray}\label{plusminus}
(T+j)-j &=& T \ \ \mbox{ if } T_{ij}=T_{i(j+1)} \ \
 \mbox{ for some } i = 1,2, \\
(T-j)+j &=& T \ \ \mbox{ if } T_{ij}<T_{i(j+1)} \ \ \mbox{ for both } i = 1,2.
\nonumber
\end{eqnarray}
Now, if $S=\{x_1<x_2< \cdots < x_t\}$, where $t=2s-r$ and $P \in \ST_{2,s}$
we let
$$
\Psi(S,P):= P-x_1-x_2-\cdots-x_t.
$$
By \eqref{plusminus} it follows that $\Psi$ is the inverse to $\Phi$ and
the theorem follows.
\end{proof}

We now have a combinatorial proof of the following theorem given
in a different form in \cite{Burstein}:
\begin{theorem}\label{closed}
For all $n,k \geq 0$ we have
$$
|[k+2]^n(123)|= \sum_{r,i}2^{n-r}C_i\binom i {r-i} \binom n r \binom k i,
$$
where $C_i$ is the $i$th Catalan number. The generating
function $$F(x,y)~:=\sum_{n,k}|[k+2]^n(123)|x^ky^n,$$ is given by
$$
F(x,y)= \frac 1 {(1-x)(1-2y)}
C\left( \frac{xy(1-y)}{(1-x)(1-2y)^2} \right),
$$
where $C(z)$ is the generating function for the Catalan numbers. Equivalently,
$F(x,y)$ is algebraic of degree two and satisfies the equation:
$$
x(1-x)y(1-y)F^2-(1-x)(1-2y)F+1=0.
$$
\end{theorem}

To complete the picture for permutation patterns of length $3$ it
remains to enumerate $|[k]^n(132)|$.
Simion and Schmidt~\cite{SimSch1985} introduced a simple
bijection between $S_n(123)$ and $S_n(132)$ which fixes each
element of $S_n(123)\cup S_n(132)$. West~\cite{Wes1990}
generalized this bijection to obtain a bijection between $S_n(p)$
and $S_n(q)$ where $p(\ell)=q(\ell-1)=\ell$,
$p(\ell-1)=q(\ell)=\ell-1$, and $p,q\in S_\ell$. Here
we indicate how to generalize West's result to obtain a bijection
between $[k]^n(p)$ and $[k]^n(q)$ where $p$ and $q$ are as above.
% Following~\cite{Wes1990} we define as follows.
% \begin{definition}
% Let $p=(p_1,\ldots,p_\ell)$ be any pattern on $d$ letters. Also
% let $\pi=(\pi_1,\ldots,\pi_n)\in [k]^n(p)$. For $1\leq j\leq n$,
% if $s$ is the largest integer for which there exists a subsequence
% $\pi_{i_1},\ldots,\pi_{i_s}=\pi_j$ which is of the same type as
% $p_1,\ldots,p_s$, then let $p_j$ be a member of the $s$th basic
% subsequence of $\pi$ with respect of $p$.
% \end{definition}
% We say that a word pattern $p$ of length $\ell$ on $d$ letters (a
% {\em word pattern} on $d$ letters is a word with repetition
% letters and every letter $1,2,\ldots,d$ appears at least once) is
% of {\em type one} if $p=(p_1,\ldots,p_{\ell-2},d-1,d)$ and of {\em
% type two} if $p=(p_1,\ldots,p_{\ell-2},d,d-1)$, where $p$ occurs
% the letter $d$ exactly once. For any word pattern $p$ of length
% $\ell$, we define a word pattern $\widetilde{p}$ by
% $\widetilde{p}=(p_1,\ldots,p_{\ell-2},p_\ell,p_{\ell-1})$.

\begin{theorem}\label{bij1221}
Let $p=p_1p_2\cdots p_\ell$ be a pattern with greatest entry equal to $d$ and
$p_{\ell-1}=d-1$, $p_\ell=d$. If $d$ occurs exactly once in $p$ then
$$
|[k]^n(p)|=|[k]^n(\widetilde{p})|,
$$
where $\widetilde{p}=p_1p_2\cdots p_{\ell}p_{\ell-1}$.
\end{theorem}
\begin{proof}
The proof is a straight forward generalization of West's algorithm presented
in \cite[Section~3.2]{Wes1990}.
\end{proof}
For example, if $p=132$ then $\widetilde{p}=123$. Hence, by
Theorem~\ref{bij1221} we get that if $p$ and $q$ are any permutation patterns
of length $3$ then $|[k]^n(p)|=|[k]^n(q)|$ for
all $n,k\geq 0$
(see~\cite{Burstein} for an analytical proof). If
$p=1232$ the $\widetilde{p}=1223$. Hence, Theorem~\ref{bij1221}
gives $|[k]^n(1232)|=|[k]^n(1223)|$ for all $n,k\geq0$.

Since, $S_n(p) \subset [n]^n(p)$, the numbers $|[n]^n(p)|$ are interesting.
A sequence $f(n)$ is {\em polynomially recursive} ({\em P-recursive}) if
there is a finite number of polynomials $P_i(n)$ such that
$$
\sum_{i=0}^NP_i(n)f(n+i)=0,
$$
for all integers $n \geq 0$.
For the case when $p$ is permutation pattern of length $3$ we have the
following:
\begin{theorem}\label{diagonal}
Let $p$ be a permutation pattern of length $3$. Then the sequence
$f(n):=|[n]^n(p)|$ is $P$-recursive and satisfies the three term
recurrence:
    $$p(n)f(n-2)+q(n)f(n-1)+r(n)f(n)=0,$$
where
\begin{eqnarray*}
p(n)&=& 3(n-3)(n-1)(3n-5)(3n-4)(5n-4),\\
q(n)&=& 288 - 1440n + 2780n^2 - 2435n^3 + 976 n^4 - 145n^5,\ \mbox{ and }\\
r(n)&=& 2(n-2)^2n(n+1)(5n-9).
\end{eqnarray*}

\end{theorem}
\begin{proof}
The fact that $f(n)$ is $P$-recursive follows easily from the
expansion of $f(n)$ as a double sum using Theorem \ref{closed} and
the theory developed in \cite{lipshitz}. The polynomials
$p,q$ and $r$ were found using the package MULTISUM (see~\cite{WR}) developed
by Wegschaider and Riese.
\end{proof}

\begin{corollary}\label{3asy}
The asymptotics of $f(n)=|[n]^n(123)|$ is given by
$$
f(n) \sim C n^{-2}\left( \frac {13} 2 \right)^n,
$$
where $C>0$ is a constant.
\end{corollary}
\begin{proof}
This is a direct consequence of Theorem \ref{diagonal} and the theory of
asymptotics for $P$-recursive sequences, see \cite{ZeilWimp}.
\end{proof}
A consequence of this is that the generating function of $f(n)$ is
transcendent, since the exponent of $n$ in the asymptotic
expansion of a sequence with an algebraic generating function is
never a negative integer.

\subsection{Generating function approach}
In this section we will investigate the generating function that
enumerates the number of segmented tableaux according to
size of rows and number of different entries. Let
$A_{\ell}(x_1,x_2, \ldots, x_\ell,t)$
be the generating function:
$$
A_{\ell} = \sum_{T} x_1^{\lambda_1(T)}x_2^{\lambda_1(T)-\lambda_2(T)}\cdots
                    x_{\ell}^{\lambda_{\ell-1}(T)-\lambda_\ell(T)}t^{N(T)},
$$
where $\lambda_i(T)$ denotes the size of row $i$ in $T$, $N(T)$ denotes the
number of different entries in $T$ and the sum is over all segmented tableaux
with at most $\ell$ rows. For $i = 1,2, \ldots, \ell$ let
$A^i_\ell(x_1,\ldots,x_\ell,t)$ be the generating function for those
tableaux which have their maximal entry in row $i$. If
$F(x_1, x_2, \ldots, x_n)$ is a formal power-series in $n$ variables
the {\em divided difference} of $F$ with respect to the variable $x_i$ is
$$
\Delta_iF:= \frac{F-F(x_i=0)}{x_i},
$$
where $F(x_i=0)$ is short for
$F(x_1,x_2,\ldots, x_{i-1},0,x_{i+1}, \ldots, x_n)$.
\begin{theorem}
With definitions as above we have that $A_\ell$ satisfies the
following system of equations:
\begin{eqnarray*}
A_\ell &=& 1 + A^1_\ell+ \cdots + A^\ell_\ell, \\
A^1_\ell &=& x_1x_2tA_\ell + x_1x_2A^1_\ell,\\
A^2_\ell &=& x_3t\Delta_2A_\ell+
             x_3\Delta_2A^2_\ell,\\
         &\vdots&\\
A^{\ell-1}_\ell &=& x_\ell t\Delta_{\ell-1}A_\ell+
     x_\ell\Delta_{\ell-1}A^{\ell-1}_\ell,\\
A^\ell_\ell &=& t\Delta_{\ell}A_\ell+
             \Delta_\ell A^\ell_\ell.
\end{eqnarray*}
\end{theorem}
\begin{proof}
The theorem follows by treating two separate cases. Let $n$ be the
greatest entry in the tableau $T$. The case when there is one $n$
in a row corresponds to the first summand and the case when there
are more than one $n$ in a row corresponds to the second summand.
\end{proof}
When $\ell=2$, $A=A_2$, the system boils down to:
\begin{equation}\label{kernel}
\left( (1-x_2^{-1})(1-\frac{x_1x_2t}{1-x_1x_2})-x_2^{-1}t\right)A =
1-x_2^{-1}(1+t)A(x_2=0).
\end{equation}
This equation can be solved using the so called {\em kernel method} as
described in \cite{banderier}. If we let
$$
x_2= \frac{ 1 + x_1(1+2t) - \sqrt{(1+x_1(1+2t))^2-4x_1(1+t)^2}} {2x_1(1+t)},
$$
then the parenthesis infront of $A$ in \eqref{kernel} cancels, and
we get:
$$
A(x_2=0)=\frac{ 1 + x_1(1+2t) - \sqrt{(1+x_1(1+2t))^2-4x_1(1+t)^2}}
{2x_1(1+t)^2}.
$$
By the interpretation of $a(\ell,k,r)$, we have that the
bi-variate generating function for $a(2,k,r)$ is
$(1+x_1)^{-1}A_2(x_1,1,t)$. From this and \eqref{sumbin} one may
derive an analytic proof of Theorem~\ref{closed}.

\section{Further results and open problems}

\subsection{Further directions}
Recall that the Stanley-Wilf Conjecture asserts that
for any permutation $\pi$ the limit
$\lim_{n \rightarrow \infty}|\sym_n(\pi)|^{1/n}$ exists and is finite. What
about the sequence $|[n]^n(\pi)|$?
\begin{problem}\label{wordstw}Let $\pi$ be a permutation.
Is there a constant $0<C<\infty$ such
$|[n]^n(\pi)| \leq C^n$ for all $n \geq 0$?
\end{problem}
Note that the answer to Problem \ref{wordstw} is no when $\pi$ is not a
permutation, since then $\sym_n(\pi) \subseteq [n]^n(\pi)$.
Again, Problem \ref{wordstw} is
equivalent to the statement that
$$
\lim_{n \rightarrow \infty}|[n]^n(\pi)|^{1/n},
$$
exists and is finite. This is because for all $m,n \geq 0$ we have
$$
|[n+m]^{n+m}(\pi)| \geq |[n]^n(\pi)|\cdot|[m]^m(\pi)|,
$$
so we may apply Fekete's Lemma on sub-additive sequences. See
\cite[Theorem 1]{Arratia} for details
(the proof extends to words word for word).
For permutations $\pi \in \sym_3$ we have by Corollary \ref{3asy} that
$\lim_{n \rightarrow \infty}|[n]^n(\pi)|^{1/n}=13/2$ as opposed to
$\lim_{n \rightarrow \infty}|\sym_n(\pi)|^{1/n}=4$.
% Another equivalent statement is
% that given a pattern $\pi$ there are constants $C,D >0$ such that
% $|[k]^n(\pi)| \leq C^kD^n$ for all $n,k \geq 0$.

For which permutations do we know that Problem \ref{wordstw} is
true? It follows from the work in \cite{AloFri2000} that  Problem
\ref{wordstw} is true for all permutations which can be written as
an increasing sequence followed by a decreasing. Also, with no
great effort B\'ona's proof \cite{Bon1999} of the Stanley-Wilf
conjecture for layered patterns may be extended to this setting.
Thus for all classes that the Stanley-Wilf conjecture is known to
hold, the seemingly stronger Problem \ref{wordstw} holds. The
following conjecture therefore seems plausible:
\begin{conjecture}
For all permutations $\pi$ we have:
$$
\exists C \forall n ( |[n]^n(p)| \leq C^n ) \Leftrightarrow
\exists D \forall n  ( |\sym_n(p)| \leq D^n ).
$$
\end{conjecture}

There are several problems concerning the automatons associated to a pattern
that has connections to the above problems. One problem is to give
an estimate to the number of simple paths in $\Au(p,k)$, another
is to estimate the number of equivalence classes in $\Au(p,k)$. Yet
another problem is to give an estimate to the maximum size of an equivalence
class.

% \begin{problem}
% Estimate the number of paths in $\Au(p,k)$ from $\langle \e \rangle$
% with no loops.
% \end{problem}
% \begin{problem}\label{pr_ek}
% Let $e_k(p)$ denote the number of states in $\Au(p,k)$. Does
% $e_k(p)$ grow at most exponentially with $k$?
% \end{problem}

% \begin{problem}
% Estimate the maximum size of an equivalence class.
% \end{problem}

\subsection{Formula for $|[k]^n(p)|$}
Our algorithm (see Theorem~\ref{gf}) for finding a formula for
$|[k]^n(p)|$ is implemented in C++ and Maple, see~\cite{M}. The first
with input $p$ and $k$ and output the automaton $\Au(p,k)$
and the second with input the automaton $\Au(p,k)$
and output the exact formula for $|[k]^n(p)|$. This algorithm
allows us to get an explicit formula for $|[k]^n(p)|$ where $p\in
S_k$ and $k\geq1$ are given. For example, an output for the
algorithm for $p\in S_4$ and $k=3,4,5,6$ is given by
Table~\ref{tab4}, where we define,
        $$[b_0,b_1,\ldots,b_d]_x=\sum_{j=0}^{d}b_j\binom{n}{j}x^{n-j}.$$

\begin{table}[h]\label{tab4}
{ \begin{tabular}{|l|l|l|}\hline
$p$ &  $k$ & $|[k]^n(p)|$\\ \hline\hline
$1234,1243$&  $3$ & $[1]_3$ \\
$1432,2143$  &$4$ & $[1,1,1,1]_3$ \\
  &$5$ & $[1,2,4,8,11,10,5]_3$ \\
  &$6$ & $[1,3,9,27,66,126,183,189,126,42]_3$ \\
  &$7$ & $[1,4,16,64,221,632,1478,2772,4074,4536,3612,1848,462]_3$\\
  &$8$ &  $[1,5,25,125,555,2103,6735,18075,40290,73770,109206,$\\
  &    & $\qquad127710,113850,72930,30030,6006]_3$\\
  &$9$ & $[1,6,36,216,1170,5508,22338,77688,230823,583410,$\\
  &    & $\qquad1247076,2235816,3322836,4025736,3880305,2867436,$\\
  &    & $\qquad1528956,525096,87516]_3$\\
  &$10$&$[1,7,49,343,2191,12313,60361,257407,953554,3064558,8527666,$\\
  &    & $\qquad20482462,42268534,74452378,110916091,137998861,$\\
  &    &
  $\qquad140882742,115068954,72390318,32978946,9699690,1385670]_3$  \\\hline\hline

$1324$&  $3$ & $[1]_3$ \\
  &$4$ & $[1,1,1,1]_3$ \\
  &$5$ & $[1,2,4,8,11,10,5,1]_3$ \\
  &$6$ & $[1,3,9,27,66,126,183,197,152,80,26,4]_3$ \\
  &$7$ & $[1,4,16,64,221,632,1478,2808,4308,5295,5152,3895,2219,$\\
  &    & $\qquad904,239,33,1]_3$\\ \hline\hline

$1342$&  $3$ & $[1]_3$ \\
  &$4$ & $[1,1,1,1]_3$ \\
  &$5$ & $[1,2,4,8,11,10,4]_3$ \\
  &$6$ & $[1,3,9,27,66,126,176,168,96,24]_3$ \\ \hline\hline

$1423$&  $3$ & $[1]_3$ \\
  &$4$ & $[1,1,1,1]_3$ \\
  &$5$ & $[1]_2+[0,3,3,9,10,11,3]_3$ \\
  &$6$ & $[13,1]_2+[-12,15,-2,37,57,134,169,167,76,12]_3$ \\  \hline\hline

$2413$&  $3$ & $[1]_3$ \\
  &$4$ & $[1,1,1,1]_3$ \\
  &$5$ & $[10,4,1]_2+[-9,8,1,9,11,10,2]_3$ \\
  &$6$ & $[96,28,5]_2+[-95,71,-36,54,52,132,167,137,44,4]_3$ \\ \hline
\end{tabular}}
\caption{Patterns of length $4$}
\end{table}
Finally we remark that our method can be generalized as follows. Given
a set of patterns $T$ we define an equivalence relation $\sim_T$
on $[k]^*$ by: $v \sim_T w$ if for all words $r \in [k]^*$ we have
$$vr\mbox{ avoids } T   \ \ \mbox{ if and only if  } \ \ wr \mbox{ avoids } T,$$
where a word $u$ avoids $T$ if $u$ avoids all patterns in $T$. As
in Section~2 we define an automaton $\Au(T,k)$ with the
equivalence classes of $\sim_T$ as states. With minor changes in
the proof, Theorem~\ref{gf} can be
extended to avoidance of a set of patterns. For example, if
$T=\{1234,2134\}$ and $k=5$, then by ~\cite{M} we get that
$$\begin{array}{l}
|[4]^n(T)|=2\cdot3^n+2\binom{n}{2}3^{n-2}-2^n,\\[4pt]
|[5]^n(T)|=3\cdot3^n+6\binom{n}{2}3^{n-2}+6\binom{n}{3}3^{n-3}+8\binom{n}{4}3^{n-4}-2\cdot2^n,\\[4pt]
|[6]^n(T)|=4\cdot3^n+12\binom{n}{2}3^{n-2}+24\binom{n}{3}3^{n-3}\\
\qquad\quad\qquad\qquad\qquad\qquad\qquad+54\binom{n}{4}3^{n-4}+60\binom{n}{5}3^{n-5}+40\binom{n}{6}3^{n-6}-3\cdot2^n.
\end{array}$$
{\bf Acknowledgements.} We would like to thank M. Bousquet-M\'elou
and C. Krattenthaler for interesting discussions, D. Zeilberger
for informing us about the work of Wegschaider and Riese, and
Wegschaider and Riese for running their program for us.
\end{document}